\documentclass[12pt]{article}

\usepackage{amsmath,amsthm,bbm}
\usepackage{version}
\usepackage{color}
\usepackage{amstext}
\usepackage{amssymb}
\usepackage{latexsym}
\usepackage{amsfonts}
\usepackage{graphicx}

\newcommand{\rank}{{\rm rank\ }}

\let\cal\mathcal

\usepackage{latexsym}
\usepackage{amssymb}
\usepackage{euscript}

\def\mcc{M\raise.5ex\hbox{c}C}
\def\mccarthy{M\raise.5ex\hbox{c}Carthy}
\def\Hu{\H}
\def\sz{Szeg\Hu{o} }

\def\eg{{\it e.g. }}
\def\ie{{\it i.e. }}

\def\l{\lambda}

\def\vare{\varepsilon}

\let\ii=\i
\let\i=\infty

\def\la{\langle}
\def\ra{\rangle}
\def\={\ = \ }

\def\J{{\cal J}}

\def\C{\mathbb C}

\def\T{\mathbb T}
\def\D{\mathbb D}
\def\B{\mathbb B}

\def\inn{\ \in \ }

\def\dis{\displaystyle}

\def\be{\setcounter{equation}{\value{theorem}} \begin{equation}}
\def\ee{\end{equation} \addtocounter{theorem}{1}}
\def\beq{\begin{eqnarray*}}
\def\eeq{\end{eqnarray*}}
\def\se{\setcounter{equation}{\value{theorem}}} 
\def\att{\addtocounter{theorem}{1}}

\def\bs{\vskip 12pt}

\def\bp{{\sc Proof: }}
\def\ep{{}{\hfill $\Box$} \vskip 5pt \par}

\def\bl{\begin{lemma}}
\def\el{\end{lemma}}
\def\bt{\begin{theorem}}
\def\et{\end{theorem}}
\def\bprop{\begin{prop}}
\def\eprop{\end{prop}}
\def\bd{\begin{definition}}
\def\ed{\end{definition}}
\def\br{\begin{remark}}
\def\er{\end{remark}}
\def\bexer{\begin{exercise}}
\def\eexer{\end{exercise}}
\def\bfig{\begin{figure}}
\def\efig{\end{figure}}

\newtheorem{theorem}{Theorem}[section]
\newtheorem{prop}[theorem]{Proposition}
\newtheorem{lemma}[theorem]{Lemma}

\newtheorem{question}[theorem]{Question}
\newtheorem{definition}[theorem]{Definition}

\theoremstyle{definition}
\newtheorem{example}[theorem]{Example}
\newtheorem{remark}[theorem]{Remark}

\newcommand{\Ci}{\C_\i}
\newcommand{\ci}{\C_\i}
\newcommand{\G}{\Gamma}
\renewcommand{\J}{{\mathbb J}}
\renewcommand{\d}{{\delta}}

\begin{document}
\setlength{\baselineskip}{21pt}

\title{The Takagi problem on the disk and bidisk}
\author{Jim Agler
\thanks{Partially supported by National Science Foundation Grant
DMS 1068830}
\and
Joseph A. Ball
\and
John E. M\raise.5ex\hbox{c}Carthy
\thanks{Partially supported by National Science Foundation Grant DMS 0966845}}

\date{August 13, 2012}

\bibliographystyle{plain}
\maketitle
\begin{abstract}
We give a new proof on the disk that a Pick problem can be solved by a rational
function that is unimodular on the unit circle and for which the number of poles inside
the disk is no more than the number of non-positive eigenvalues of the Pick matrix. 
We use this method to find rational solutions to Pick problems on the bidisk.
\end{abstract}

{\em Dedicated to the memory of Bela Szokefalvi-Nagy, whose work inspired each of us.}

\baselineskip = 18pt

\section{Introduction}\label{seca}

Given points $(\l_1,\dots,\l_N)$ in the unit disk $\D$, and complex numbers $(w_1,\dots,w_N)$,
a classical theorem of G.~Pick \cite{pi16} asserts that there exists a holomorphic function
$\phi$ on $\D$ that interpolates the data, \ie 
satisfies
\be
\label{eqxa1}
\phi(\l_i) \= w_i \qquad \forall\ 1 \leq i \leq N
\ee
and, in addition, maps $\D$ to $\D$, if and only if the Pick matrix
\be
\label{eqa1}
\G \= \left( \frac{1-w_i \bar w_j}{1- \l_i \bar\l_j} \right)_{i,j=1}^N
\ee
is positive semi-definite. Moreover, when $\G$ is positive semi-definite, $\phi$ can be chosen to 
be a Blaschke product with degree equal to the rank of $\G$, and so extends to be meromorphic on the whole Riemann
sphere $\Ci$, and to have modulus one on the unit circle $\T$.

The case where $\G$ has some negative eigenvalues was first studied by
T.~Takagi \cite{tak29}, and later by many other authors
\cite{aak71,nud77,bh83,bkr03}. See the book \cite{bgr90} for an account.
The principal difference is that if one wishes to interpolate with a unimodular function
(\ie a function that has modulus one on $\T$), then one has to allow poles inside $\D$. 
For a rational function $\phi = q/p$, let us say that $\phi$ {\em strictly satisfies the interpolation condition at $\l_i$} if
\be
\label{eqxa2}
\lim_{\l \to \l_i} \phi(\l) \= w_i ;
\ee
and we shall say that 
$\phi$ {\em weakly satisfies the interpolation condition at $\l_i$} if
\be
\label{eqxa3}
 q(\l_i) \= w_i p(\l_i) .
\ee
If $p(\l_i) \neq 0$, clearly these two conditions are the same.

A typical result, see \eg \cite[Theorem 19.2.1]{bgr90} ,
is
\bt
\label{thmaa}
Suppose the matrix $\G$ in (\ref{eqa1}) is invertible, and has $\pi$ positive eigenvalues and
$\nu$ negative eigenvalues. Then there exists a meromorphic 
function $\phi$ that satisfies the strong interpolation conditions,
is unimodular, and is the quotient of a Blaschke product of degree
$\pi$ by a Blaschke product of degree $\nu$.
\et

The identity 
\be
\label{eqa19}
1 - f/g \= (1/g) \left[ (1-f) - (1-g) \right]
\ee
shows that one cannot use Blaschke products of lower degree, as
this would result in  $\Gamma$ having at most deg $f$ positive eigenvalues,
and at most deg $g$ negative eigenvalues.

If $\G$ is not invertible, the problem is degenerate, and can be treated as a limiting case
of non-degenerate problems. One way to do this is to seek a unimodular
function with only $\nu$ poles, but which may not have the value $w_i$ at every $\l_i$ if
both the numerator and denominator had zeroes there before cancelling out those factors.
See the paper \cite{bkr03} for recent developments on this approach.

The strict interpolation problem was solved by H. Woracek \cite{wor95};
see also V. Bolotnikov's paper \cite{bol09}.

\bt
\label{thmab}
Suppose 
$\l_1,\dots,\l_N\inn \D$ and $w_1,\dots,w_N \inn \C$.
Suppose the matrix
\be
\label{eqg1}
\G \= \left( \frac{1-w_i \bar w_j}{1- \l_i \bar\l_j} \right)_{i,j=1}^N
\ee
has $\pi$ positive eigenvalues, $\nu$ negative eigenvalues, and $\zeta$ zero eigenvalues.
Then there exists a rational function $\phi$ that is unimodular on $\T$,
such that $\phi(\l_i) = w_i$ for every $i$, and such that 
$\phi$ has at most $N- \pi$ poles and $N-\nu$ zeroes in $\D$.

\et

%

In Section~\ref{secb}, we 
prove Theorem~\ref{thmab} using the idea of lurking isometries.
Our proof depends on the idea of considering an indefinite inner product that
is associated with the problem.
This idea was first introduced by M.G.~Kre{\u{\ii}}n \cite{kr70}, and
was further developed in the context of the Takagi problem 
in the Grassmannian approach of
 J.A.~Ball and J.W.~Helton \cite{bh83}.

In Section~\ref{secc} we prove a bidisk version of Theorem~\ref{thmab},
though with weaker bounds on the degrees.

\section{Degenerate Case on Disk}
\label{secb}
If $p$ is a polynomial of degree $d$, let 
\[
\tilde p(z) \= z^d \overline{p(1/\bar z)} .
\]

Proof of Theorem~\ref{thmab}:

{\bf Step 1:} Assume $\l_1 = 0$.

There are vectors $u_i$ in $\C^\pi$ and $v_i$ in $\C^\nu$ such that
\be
\label{eqb1}
\G_{ij} \= \la u_i, u_j \ra \ - \ \la v_i, v_j \ra .
\ee
Rewrite (\ref{eqb1}) 
as
\be
\label{eqb2}
1 - w_i \bar w_j \= (1 - \l_i \bar \l_j) \left[ \la u_i, u_j \ra \ - \ \la v_i, v_j \ra \right].
\ee
 
Choose vectors $y_i$ in $\C^{N-\pi - \nu}$ so that the matrix
$$
B_{ij} \ := \ \la u_i, u_j \ra \ + \ \la v_i, v_j \ra \ + \ \la y_i, y_j \ra
$$
is positive definite.

Let
$$
x_i \= \left( \begin{array}{c}
u_i \oplus y_i \\
v_i \oplus y_i
\end{array} \right) 
$$
be vectors in $\C^{2N-\pi - \nu}$. 
Write $N' = 2N - \pi - \nu$.
Then the vectors $\{ x_i \}_{i=1}^N$ are linearly independent vectors
in $\C^{N'}$, since their Gramian has rank $N$.

Let $J$ be the signature matrix
$$
J \= 
\bordermatrix{&\C  &\C^{N-\nu}& \C^{N-\pi} \cr
\C &1 & 0&0\cr
\C^{N-\nu}&0&I &0\cr
\C^{N-\pi}&0&0&-I} ,
$$
and $J_1$ be the compression of $J$
$$
J_1 \= 
\bordermatrix{&\C^{N-\nu}& \C^{N-\pi} \cr
\C^{N-\nu}&I &0\cr
\C^{N-\pi}&0&-I} .
$$

Then (\ref{eqb2}) can be written as
\be
\label{eqb3}
\la J 
\left( \begin{array}{c} 1 \\  \l_i x_i \end{array} \right) ,
\left( \begin{array}{c} 1 \\  \l_j x_j \end{array} \right) \ra
\=
\la J 
\left( \begin{array}{c}  w_i \\  x_i \end{array} \right) ,
\left( \begin{array}{c}  w_j \\  x_j \end{array} \right) \ra .
\ee

Define the linear map $V$ by
$$
V\ : \ \sum_{i=1}^N c_i \left( \begin{array}{c} 1 \\  \l_i x_i \end{array} \right)
\ \mapsto \
\sum_{i=1}^N c_i \left( \begin{array}{c}  w_i \\  x_i \end{array} \right) .
$$
Since the $x_i$'s are linearly independent, 
the map $V$ is a well-defined linear map from an 
$N$ dimensional subspace of $\C^{N' + 1}$
onto another 
$N$ dimensional subspace of $\C^{N' + 1}$.
Moreover, from (\ref{eqb3}), we have that 
$$
V^\ast J V \= J .
$$
Because $V$ is one-to-one, it can be extended to a map $V_1 : \C^{N'+1} \to
\C^{N' + 1}$ that is still a $J$-isometry, \ie

\be
\label{eqb355}
V_1^\ast J V_1 \= J.
\ee
(See \eg \cite[p.264]{azi89}).
Write
\be
\label{eqb37}
V_1 \= 
\bordermatrix{&\C  &\C^{N'} \cr
\C &A & B\cr
\C^{N'}&C&D }
\ee
and, since $V_1$ extends $V$,
\be
V_1
\ : \ 
 \left( \begin{array}{c} 1 \\ \l_i x_i \end{array} \right)
\ \mapsto \
 \left( \begin{array}{c}  w_i \\  x_i \end{array} \right) .
\label{eqb4}
\ee
Define $\phi$ by
\be
\label{eqb47}
{\phi(\l)} \=
A +  \l B (I -  \l D)^{-1} C .
\ee
From (\ref{eqb4}), we get
\be
\label{eqxb1}
(I -  \l_i D) x_i \= C 1 ,
\ee
so if $(I -  \l_i D)$ is invertible, we get that $\phi(\l_i) = w_i $;
otherwise (\ref{eqxa3}) holds.

Indeed, in the latter case, write $\phi$ as
\[
\phi(\l) \= \frac{1}{\det(I - \l D)}
\left[ \det(I - \l D) A + \l B \Xi(\l) C \right]  \ =: \  \frac{q}{p}(\l),
\]
where $\Xi(\l)$ is the adjugate matrix  of $(I - \l D) $.
If $(I - \l_i D)$ is not invertible, then $p(\l_i) = 0$; but using (\ref{eqb37}) and (\ref{eqxb1}), we get
\be
\label{eqxb5}
\phi(\l)
= w_i +  (\l - \l_i) B (1-\l D)^{-1}Dx_i ,
\ee
so the order of the pole of $\phi$ at $\l_i$ is 
one less than the order of the zero of $p$, so $q$ must also vanish at $\l_i$.

Note that at $\l_1$ we get the strict interpolation condition, because
we have assumed $\l_1 = 0$ so $(I -  \l_1 D)$ is invertible.

From (\ref{eqb47}), we get
$$
1 - |\phi(\l)|^2 \= (1-|\l|^2) C^* (I -  \bar \l D^*)^{-1} J_1  (I -  \l D)^{-1} C .
$$
Therefore 
$$
\lim_{|\l| \nearrow 1} \, 1 - |\phi(\l)|^2 \= 0 ,
$$
except possibly at the finite set $\sigma(D) \cap \T$. As $\phi$ is rational, it follows that it must
be unimodular, and therefore a ratio $\tilde p/p$ of a polynomial with its reflection,  where
$$
{\rm deg}(p)  \ \leq \ N' ,
$$
and $p$ has no zeroes on the unit circle $\T$.
We can write 
\[
\tilde p/p \= f/g ,
\]
where $f$ and $g$ are Blaschke products whose zeroes correspond to those of $p$ in ${\mathbb E} = 
\C \setminus \overline{\D}$ and $\D$ respectively.

For $\l$ not in $1/\sigma(D)$, define
\be
\label{eqb7}
x(\l) \= (I -  \l D)^{-1} C 1 .
\ee
Then
\be
\label{eqb75}
\bordermatrix{&  & \cr
 &A & B\cr
&C&D }
\  
 \left( \begin{array}{c} 1 \\  \l x(\l) \end{array} \right)
\=
 \left( \begin{array}{c}  \phi(\l) \\  x(\l) \end{array} \right) .
\ee
Let $\mu_1,\dots,\mu_k$ be arbitrary points in $\D \setminus 1/\sigma(D)$.
By (\ref{eqb75}), 
\be
\label{eqb8}
\la J 
 \left( \begin{array}{c} 1 \\  \mu_i x(\mu_i) \end{array} \right)
,
 \left( \begin{array}{c} 1 \\  \mu_j x(\mu_j) \end{array} \right) \ra
\=
\la J
 \left( \begin{array}{c}  \phi(\mu_i) \\  x(\mu_i) \end{array} \right) ,
 \left( \begin{array}{c} \phi(\mu_j) \\  x(\mu_j) \end{array} \right) \ra.
\ee
Rewriting (\ref{eqb8}), we get
\be
\label{eqb9}
\frac{1 - \phi(\mu_i) \bar \phi(\mu_j)}{1- \mu_i \bar \mu_j}
\= 
\la J_1 x(\mu_i), x(\mu_j) \ra .
\ee
So a maximal negative eigenspace of
$$
\G(\l,\mu) \= 
\frac{1 - \phi(\l) \bar \phi(\mu)}{1- \l \bar \mu}
$$
can have no larger dimension than that of $J_1$, which is $N-\pi$, and
a maximal positive eigenspace 
can have dimension at most $N-\nu$.

However, if we choose the $\mu_i$'s to be the set $\{ \mu \inn \D \ : \
\phi(\mu) = M \}$ for any $|M| > 1$, we get a negative definite subspace for $\Gamma$ of dimension
${\rm deg}(g)$, and if we choose $|M| < 1$, we get a positive definite subspace of dimension
${\rm deg}(f)$. Therefore 
${\rm deg}(f) \leq N - \nu$, and  
${\rm deg}(g) \leq N - \pi$.

One could also argue that as $D$ is a $J_1$-contraction, the number of points in
$\sigma (D) \cap \overline{\D}$ is at most $\nu + \zeta$, and the number
in $\sigma (D) \cap \overline{\mathbb E}$ is at most $\pi + \zeta$;
this follows from  \cite[Thm. 4.6.1]{glr05}
and a perturbation argument.

{\bf Step 2:}
Let $m_j$ be the Mobius map that swaps $0$ and $\l_j$.
Applying step 1 to the modified problem
\[
m_j(\l_i) \to w_i, \qquad 1 \leq j \leq N
\]
we get a solution as in Step 1, and postcomposing this solution with $m_j$ we get functions
$\phi_j$ such that:

$\bullet$ Each $\phi_j$ solves the weak interpolation problem at every point, and satisfies  the strong interpolation condition at $\l_j$.

$\bullet$ Each $\phi_j$ is a ratio of Blaschke products of degrees at most $N-\nu$ and $N-\pi$.

Each $\phi_j$ can be written as
\[
\phi_j \= \frac{\tilde p_j}{p_j} ,
\]
where $p_j$ is a polynomial of degree less than or equal to $N'$, and with at most
$N-\nu$ zeroes in $\D$ and $N-\pi$ zeroes in ${\mathbb E} $.
Multiplying $p$ (and hence $\tilde p$) by an appropriate power of $(1+z)$, we can assume that
every $p_j$ has {\em exactly} the same degree, some number $ \leq N'$.

{\bf Step 3:}
Let $t_j$ be real numbers such that
\[
\sum t_j p_j (\l_i) \ \neq \ 0, \qquad \forall \ 1 \leq i \leq N .
\]
(They exist since $p_i(\l_i) \neq 0$.)
Let $q = \sum t_j p_j$.
Then $\phi = {\tilde q}/q$ is a ratio of Blaschke products, and satisfies the strict interpolation 
condition at every node.

The number of zeroes plus the number of poles of $\phi$ in $\D$ is at most $N'$, and
as the Pick matrix of $\phi $ is $\Gamma$, it has at least $\pi$ zeroes and at least $\nu$ poles.

 Write $\phi = f/g$ where $f$ and $g$ are Blaschke products with
no common factors.
We have $\deg f \geq \pi$ and $\deg g \geq \nu$, and we wish to prove that
$\deg f \leq \pi + \zeta$ and $\deg g \leq \nu + \zeta$. If 
$N'' := \deg f + \deg g \leq N = \pi + \nu + \zeta$, then
\[
(\deg f - \pi) + (\deg g - \nu) \leq \zeta.
\]
As we know that each term on the left-hand side is non-negative, we get that both
\[
\deg f - \pi \leq \zeta, \qquad \deg g - \nu \leq \zeta,
\]
and hence
\[
\deg f \leq \pi + \zeta, \qquad \deg g \leq \nu + \zeta
\]
and we are done.

So we shall assume that $N'' > N$.

{\bf Step 4}: Choose points $\l_{N+1}, \dots, \l_{N''}$ in $\D$ that are
not in the zero set of $g$,
 and so that the Pick matrix for $\phi$ at the nodes
$\{\l_{N+1}, \dots, \l_{N''} \}$ is invertible.
We claim that this can be done by choosing points in the level set $\{ \phi = c \}$
for any constant $c$ which is not unimodular. Indeed, if $|c| < 1 $, a homotopy argument based
on the argument principle (and using that $\phi$ is unimodular on the unit circle) shows that the
number of points of $\D$ 
in the level set $\{ \phi = c \}$ is the same as the number of zeroes of $\phi$ in $\D$, namely
$\deg f$. Similarly, if $|c| > 1$, then the number of points 
 of $\D$ 
in the level set $\{ \phi = c \}$ is the same as the number of poles of $\phi$ in $\D$, namely
$\deg g$. If $\deg f < N'' - N$, then $\deg f < \deg f + \deg g - (\pi + \nu + \zeta)$,
implying that $\deg g > \zeta = N' - N \geq N''- N $. Similarly, $\deg g < N'' - N$ forces
$\deg f > N'' - N$. We conclude that at least one of $\deg f$ and $\deg g$ is at least $N'' - N$.
It follows that there is a choice of non-unimodular $c$ so that the level set $\{ \phi = c\}$
has $N'' - N$ points $\lambda_{N+1}, \dots, \lambda_{N''}$ as required.


By Lemma~\ref{lemb1} below, the Pick matrix $\Pi$ for $\{\l_1, \dots, \l_{N''} \}$ is invertible,
and has $\deg f$ positive eigenvalues and $\deg g$ negative ones. We can write it in block
form as
\[
\Pi \=
\left(
\begin{matrix}
\Gamma&B\\B^*&C
\end{matrix}
\right)
\]
where $C$ is invertible.
The inertia of $\Pi$ (the numbers of positive, negative and zero eigenvalues)
then equals the inertia of $C$ plus the inertia of the Schur complement,
$\Gamma - B C^{-1} B^*$.
Moreover, the size of $C$ is $\eta$-by-$\eta$, where $\eta = N''-N \leq  \zeta$.
If $C$ has inertia $(\eta_1,\eta_2, 0)$, where $\eta_1 + \eta_2 = \eta$,
then $ - B C^{-1} B^*$ has inertia $(\eta_2,\eta_1,0)$.
As $\Gamma$ has a $\zeta$-dimensional null space, and $\Gamma - B C^{-1} B^*$
is invertible,
$B C^{-1} B^*$ must have rank at least $\zeta$. On the other hand,
\[
\rank B C^{-1} B^* \leq \rank C^{-1} = \eta \leq \zeta
\]
and we conclude that $\eta = \zeta$. If $C$ has inertia
$(\eta_1,\eta_2,0)$, then $\eta_1 + \eta_2 = \eta = \zeta.$
As $B C^{-1} B^*$ has rank $\zeta$, it follows that the inertia of
$-B C^{-1} B^*$ must be $(\eta_2,\eta_1, N- \zeta)$. 
As $\Gamma = B C^{-1} B^*$ is invertible, then necessarily 
$\Gamma - B C^{-1} B^*$ must have inertia $(\pi+\eta_2,
\nu+\eta_1,0)$. 
Hence $\Pi$ has inertia $(\pi+\zeta,\nu+\zeta,0)$, and the degrees of $f$ and $g$
are exactly $\pi + \zeta$ and $\nu + \zeta$ respectively.
\ep

Let $k_\l(z) = \frac{1}{1-\overline{ \l }z}$ be the \sz  kernel at $\l$ in the Hardy space $H^2$.
The following lemma was proved in \cite[Lemma 3.3]{bokh06}; we include a proof for completeness.
\bl
\label{lemb1}
Let $f$ and $g$ be relatively prime Blaschke products of degrees $m$ and $n$ respectively.
Let $\Lambda$ be any set of $m+n$ distinct points in $\D$ that is disjoint from the zero set of
$g$. Then the Pick matrix for $\phi = f/g$ at the points of $\Lambda$ has inertia $(m,n,0)$.
\el 
\bp
Let $N = m+n$, and let $\Lambda = \{ \l_1, \dots, \l_N \}$.
Let
\[
\Delta_{ij} \= \frac{g(\l_i) \overline{g(\l_j)} - f(\l_i) \overline{f(\l_j)}}{1  - \l_i \overline{\l_j}}
\= \left[ g(\l_i) \overline{g(\l_j)} \right]
 \frac{1 - \phi(\l_i) \overline{\phi(\l_j)}}{1  - \l_i \overline{\l_j}}
\]
Then $\Delta$ will have the same inertia as the Pick matrix for $\phi$.
Suppose $\Delta$ has a null-vector $(c_1, \dots, c_N )^t$. Let
$ \dis \psi (z) = \sum c_i k_{\l_i} $
and
\be
\label{eqw3}
\vartheta := (T_g T_g^* - T_f T_f^* )\psi,
\ee
where
$T_f$ is multiplication by $f$ on the Hardy space $H^2$.
Then $\vartheta$ must vanish at each point of $\Lambda$.
If $h \in H^2$, then (\ref{eqw3}) yields that
\beq
\la \vartheta, fg h \ra &\=&
\la ((T_g T_g^* - T_f T_f^* )\psi, fgh \ra \\
&=& \la T_g^* \psi, fh \ra - \la T_f^* \psi, gh \ra \\
&=& \la ( T_f^* T_g^* - T_g^* T_f^*) \psi, h \ra \\
&=& 0 .
\eeq
So
$\vartheta$ is in $fg{H^2}^\perp$,
and hence  is a linear combination of the \sz kernel functions at the $N$ zeroes of $fg$.
Therefore $\vartheta$ must be a rational function whose numerator is a polynomial
of degree at most $N-1$ (since each of the kernel functions vanishes at infinity).
If this vanishes at $N$ distinct points, it must be identically zero.

Writing $\psi = h_1 + f m_1 = h_2 + g m_2$, where $h_1 \perp fH^2$
and $h_2 \perp g H^2$ then if $\vartheta = 0$,  (\ref{eqw3}) says
$fm_1 = gm_2$, so $h_1 =h_2 = 0$, and $\psi \in fg h^2$.
But $\psi$ is a linear combination of $N$ kernel functions, so is a rational 
function whose numerator is of degree at most $N-1$; if it vanishes at the $N$ zeroes
of $fg$, it must be identically zero.

Therefore we can conclude that $\Delta$ is non-singular.
But we can write $\Delta$ as
\[
\Delta_{ij} \= \left[  \frac{1 - f(\l_i) \overline{f(\l_j)}}{1  - \l_i \overline{\l_j}}\right]
-
\left[  \frac{1 - g(\l_i) \overline{g(\l_j)}}{1  - \l_i \overline{\l_j}} \right],
\]
which is the difference of two positive semi-definite matrices of ranks $m$ and $n$.
Therefore $\Delta$ has exactly $m$ positive and $n$ negative eigenvalues.
\ep

\begin{example}
One may need to choose $f$ and $g$ to have the maximum degrees, $\pi+\zeta$ and $\nu + \zeta$ respectively.
Indeed, for $N \geq 3$,
let $\l_1,\dots,\l_N$ be distinct points in $\D$, with $\l_1 = 0$. Let $w_1 = 0$, and
$w_2 = \dots = w_N = 1$. Then the Pick matrix has inertia $(\pi,\nu,\zeta) \= (1,1,N-2)$.
If $\phi = f/g$ is a ratio of Blaschke products that interpolates, then
by the Schwarz reflection principle, $\phi$ takes the value $1$ at
$\{\l_2,\dots,\l_N,\l_2^{-1},\dots,\l_N^{-1}\}$. So $\phi$ must have degree
at least $2N-2$.

By the open mapping theorem, there exists $\vare > 0$ such that $\phi$ attains the value $1+\vare$
at least $N-1$ times in $\D$. As $\phi$ is unimodular on $\T$, it follows from the argument principle that
$\phi$ must have at least $N-1$ poles in $\D$. Likewise, it must have at least $N-1$ poles in $\ci \setminus \overline{\D}$,
so by reflection
$\phi$ must have at least $N-1$ zeroes in $\D$. Therefore both $f$ and $g$ must each have degree at least
$N-1$.
However, to find Blaschke products $f$ and $g$ that satisfy
$$
f(\l_i) \= w_i g(\l_i) ,
$$
we can take $f(z) = g(z) = z$.
\end{example}

\section{The bidisk}
\label{secc}

{\em Set-up and Notation.}
Suppose 
$(\l_1,\dots,\l_N)$ is an $N$-tuple in $\D^2$, and $(w_1,\dots,w_N)$ is an
 $N$-tuple in $\ci$.
Suppose that $\G^1$ and $\G^2$ are self-adjoint matrices satisfying
\be
\label{eqf2}
(1 - w_i \bar w_j)_{i,j=1}^N \= 
\sum_{r=1}^2 \ (1- \l_i^r \bar \l_j^r) \G^r_{ij}. 
\ee

When we write an equation with superscript $r$, we shall mean ``for $r = 1$ and $2$''.
Let
\be
\label{eqc1}
\Gamma^r_{ij} \= \la u^r_i , u^r_j \ra_{\C^{\pi^r}} 
\ - \ \la v^r_i , v^r_j \ra_{\C^{\nu^r}} .
\ee
Let 
\beq
\Lambda^r_{ij} &\=& \l^r_i \bar \l^r_j \\
W_{ij} &=& w_i \bar w_j .
\eeq
Let $\kappa^r \= \pi^r + \nu^r$, 
and let 
$\pi \= \pi^1 + \pi^2,\
\nu \= \nu^1 + \nu^2,\
\kappa \= \kappa^1 + \kappa^2$.
Let $\Delta^r$ be the rank $\kappa^r$ positive matrix
$$
\Delta^r_{ij} \= \la u^r_i , u^r_j \ra_{\C^{\pi^r}} 
\ + \ \la v^r_i , v^r_j \ra_{\C^{\nu^r}} ,
$$
and let $\cdot$ denote Schur multiplication of matrices (\ie entrywise multiplication).
Let $\J$ denote the matrix all of whose entries are $1$.

We shall distinguish between two cases, which are analogous to whether
or not
the matrix in (\ref{eqg1}) is non-singular.

Case (1): Suppose

(a) The rank of $W + \Delta^1 + \Delta^2$ is $N$;

(b) The rank of $\J + \Lambda^1 \cdot \Delta^1 + \Lambda^2 \cdot \Delta^2$ is $N$.

\bs

Case (2): There are positive semi-definite $N$-by-$N$ matrices $Y^r$ of rank $\d^r$ 
(and we shall write $\d = \d^1 + \d^2$)
such that

(a) The rank of $W + \Delta^1 + Y^1 + \Delta^2 + Y^2 $ is $N$;

(b) The rank of $\J + \Lambda^1 \cdot (\Delta^1  + Y^1) + \Lambda^2 \cdot (\Delta^2 + Y^2)$ is $N$.

\bs
A polynomial is called
a {\em toral polynomial} if the intersection of its zero set with $\T^2$ 
is finite. A {\em balanced disk}
is a one-dimensional analytic subvariety of $\D^2$
of the form
$\{ (z,m(z)) \ : \ z \inn \D \}$ for some M\"obius map $m: \D \to \D$.
See \cite{ams06} for more information on toral polynomials, and
\cite{agmc_vn} 
for more information on balanced disks.

\bt
\label{thmc5}
Let the notation be as above.
Then there is a rational function $\phi$ of bidegree at most $(\pi^1+\nu^1 + \d^1, \pi^2+\nu^2 + \d^2)$
that solves the weak interpolation
problem  and  is unimodular on $\T^2$ except for at most a finite set.
If $\phi$ is written as the ratio $p/q$ of polynomials with no common factor, then
both $p$ and $q$ are atoral.

Moreover:

(i) The matrices $\G^r$ can be extended to forms on $\D^2$ with at most
$\pi^r + \d^r$ positive and $\nu^r + \d^r$ negative eigenvalues, respectively, such that
\be
\label{eqf7}
1 - \phi(\l) \overline{\phi(\mu)} \= \sum_{r=1}^2 (1 - \l^r \bar \mu^r) \G^r(\l,\mu) .
\ee
(ii) The restriction of $\phi$ to any analytic disk has at most $\nu+\delta$
poles and $\pi+\d$ zeroes.

\et


\bp

Case (1). 
%
%
Let $\{x_i\}_{i=1}^N$ be the vectors in $\C^\kappa$ given by
$$
x_i \=
\left( \begin{array}{c}
u_i^1 \\
v_i^1 \\
u_i^2 \\
v_i^2 
\end{array} \right) .
$$
Let
$$
J_1 \= 
\bordermatrix{&\C^{\pi^1}& \C^{\nu^1} &\C^{\pi^2}&\C^{\nu^2} \cr
\C^{\pi^1}&I &0&0&0\cr
\C^{\nu^1}&0&-I&0&0\cr
\C^{\pi^2}&0 &0&I&0\cr
\C^{\nu^2}&0&0&0&-I},
$$
and
$$
J \= 
\bordermatrix{&\C& \C^{\kappa} \cr
\C&1 &0\cr
\C^{\kappa}&0&J_1} .
$$
Finally, for $\l \= (\l^1,\l^2)$ let $E_\l$ be the diagonal operator
$$
E_\l \=
\bordermatrix{&\C^{\kappa^1}& \C^{\kappa^2} \cr
\C^{\kappa^1}&\l^1 I &0\cr
\C^{\kappa^2}&0&\l^2 I} .
$$

Then the equation
\be
\label{eqc2}
(1 - w_i \bar w_j)_{i,j=1}^N \= 
(1- \l_i^1 \bar \l_j^1) \G^1_{ij} \ + \
(1- \l_i^2 \bar \l_j^2) \G^2_{ij} 
\ee
can be rewritten as
\be
\label{eqc3}
\la J 
\left( \begin{array}{c} 1 \\ E_{ \l_i} x_i \end{array} \right) ,
\left( \begin{array}{c} 1 \\ E_{ \l_j} x_j \end{array} \right) \ra
\=
\la J 
\left( \begin{array}{c}  w_i \\  x_i \end{array} \right) ,
\left( \begin{array}{c}  w_j \\  x_j \end{array} \right) \ra .
\ee
Define
\be
\label{eqc4}
V\ : \ \sum_{i=1}^N c_i \left( \begin{array}{c} 1 \\ E_{ \l_i} x_i \end{array} \right)
\ \mapsto \
\sum_{i=1}^N c_i \left( \begin{array}{c}  w_i \\  x_i \end{array} \right) .
\ee
By assumption (b), the vectors
$$
\left\{ 
\left( \begin{array}{c} 1 \\ E_{ \l_i} x_i \end{array} \right) 
\right\}_{i=1}^N
$$
are linearly independent, and by assumption (a) the vectors
$$
\left\{ 
\left( \begin{array}{c}  w_i \\  x_i \end{array} \right) 
\right\}_{i=1}^N
$$
are.
Therefore, $V$ is a well-defined injective linear operator from 
an $N$-dimensional subspace
of $\C^{1+\kappa}$ onto
an $N$-dimensional subspace
of $\C^{1+\kappa}$.
Moreover, $V^* J V = J$, by (\ref{eqc3}).
Therefore $V$ can be extended to a $J$-isometry $V_1 \, :\,  
\C^{1+\kappa} \to \C^{1+\kappa} $.
Write
\be
\label{eqc5}
V_1 \= 
\bordermatrix{&\C  &\C^{\kappa} \cr
\C &A & B\cr
\C^{\kappa}&C&D }
\ee

Define $\phi$ by
\be
{\phi(\l)} \=
A + B E_{ \l} (I - D E_{ \l} )^{-1} C .
\label{eqc6}
\ee

Observe: From (\ref{eqc6}), we see that $\phi$ is a rational function of degree at most $\kappa^1$ in
the first variable and $\kappa^2$ in the second.

From (\ref{eqc4}), we get 
$\phi(\l_i) = w_i$, provided
 $ (I - D E_{ \l_i} )$ is  invertible.
Otherwise,
\[
(I - D E_{\l_i}) x_i= C 1,
\]
so
\be
\label{eqcx2}
\phi(\l) \=
A + BE_\l x_i + B E_\l (I-DE_\l)^{-1} D (E_\l - E_{\l_i}) x_i .
\ee
Writing $\phi$ as $p/q$ where $q(\l) = \det(I - DE_\l)$, we therefore have that either
$\phi(\l_i) = w_i$ or
$q(\l_i) = 0 = p(\l_i)$, as in (\ref{eqxb5}).

From (\ref{eqc6}) and the fact that $V_1^* J V_1 = J$, we get that
\be
\label{eqc7}
1 - \phi(\l) \overline{\phi(\mu)} 
\=
C^* (I - E_{\bar \mu} D^*)^{-1} \left[J_1-  E_{ \bar \mu} J_1 E_\l  \right] (I - D E_{\l} )^{-1} C .
\ee

The zero set of $q$, denoted $Z_q$, is 
the algebraic set
\be
\label{eqc8}
S \= \{ \l \inn \C^2 \ : \ \det(I - E_\l D) \= 0 \} .
\ee
On $\T^2 \setminus S$, it follows from (\ref{eqc7}) that $\phi$ is unimodular. By continuity, therefore,
\be
\label{eqc9}
\T^2 \cap Z_q \ \subset  \ Z_p .
\ee
Therefore $q$ is a {\em toral polynomial}, \ie $\T^2 \cap Z_q$ is $0$-dimensional. 
Likewise
$$
Z_p \cap \T^2 \ \subseteq \ S,
$$
and again by continuity 
$$
Z_p \cap \T^2 \ \subseteq \ Z_q ,
$$
so $p$ is also toral.

Invoking continuity again, we see that $\phi$ is unimodular on $\T^2$ except on the finite
(possibly empty) singular set 
$$
Z_p \cap Z_q \cap \T^2 
$$
where it is not defined.

\bs
Proof of (i):
Define $\G^r(\l,\mu)$ by
\beq
\G^1(\l,\mu)
&\=&
C^* (I - E_{\bar \mu} D^*)^{-1} \left(  
\begin{array}{cccc}
I&0&0&0\\
0&-I&0&0\\
0&0&0&0\\
0&0&0&0
\end{array}
\right) 
(I - D E_{\l} )^{-1} C \\
&&\\
\G^2(\l,\mu)
&\=&
C^* (I - E_{\bar \mu} D^*)^{-1} \left(  
\begin{array}{cccc}
0&0&0&0\\
0&0&0&0\\
0&0&I&0\\
0&0&0&-I
\end{array}
\right) 
(I - D E_{\l} )^{-1} C.
\eeq
Then (\ref{eqf7}) follows from (\ref{eqc7}).

\bs
Proof of (ii): 
Let $\psi(z) = \phi(z,m(z))$, where $m$ is a M\"obius map, and
let
$$
{\cal D} \= \{ (z,m(z)) \ : \ z \inn \D \} .
$$
Then from (\ref{eqf7}),
\se\att
\begin{eqnarray}
\nonumber
1 - \psi(z) \bar \psi(w)
\=
&
(1 - z \bar w) \left[
\G^1((z,m(z)),(w,m(w)) \ + \right. \\
&\quad \left. \frac{1 - m(z) \bar m(w)}{1 - z \bar w} \
\G^2((z,m(z)),(w,m(w)) \right].
\label{eqc13}
\end{eqnarray}

As long as $S \cap {\cal D}$ is finite, (\ref{eqc13}) yields that off this finite set,
the Pick matrix has at most $\pi$ positive and $\nu$ negative eigenvalues.
So repeating the argument in Step 1 of the proof of Theorem~\ref{thmab},
we get that $\psi$ has at most $\pi$ zeroes and at most $\nu$ poles.

It remains to prove that $S \cap {\cal D}$ is finite.
Write
$$
D \= 
\bordermatrix{&\C^{\kappa^1}& \C^{\kappa^2} \cr
\C^{\kappa^1}&D_{11} &D_{12}\cr
\C^{\kappa^2}&D_{21}&D_{22}} .
$$
Then the defining equation for $S$, from (\ref{eqc8}), restricted to ${\cal D}$, becomes
\be
\label{eqc14}
\det \left(
\begin{array}{cc}
I - z D_{11} & - zD_{12} \\
-m(z) D_{21} & I - m(z) D_{22}
\end{array} \right) \= 0 .
\ee
We must show that the right-hand side of (\ref{eqc14}) cannot be identically zero for any
M\"obius map $m$.
This follows from the identity
\beq
&&\left(
\begin{array}{lr}
I - z D_{11} & - zD_{12} \\
-m(z) D_{21} & I - m(z) D_{22}
\end{array} \right) \= \\
&&\qquad
\left(
\begin{array}{lr}
I  & - zD_{12}(I - m(z) D_{22})^{-1} \\
0 & I 
\end{array} \right)  \\
&&\qquad\qquad\left(
\begin{array}{lr}
I-zD_{11}   - zm(z) D_{12}(I - m(z) D_{22})^{-1}D_{21}&0 \\
-m(z) D_{21} & I - m(z) D_{22} 
\end{array} \right) 
.
\eeq

Case (2): Proceed as above, but first direct sum appropriate vectors $y^r_i$ to both
$u^r_i$ and $v^r_i$.
\ep

Repeating steps 2 and 3 of the proof of Theorem~\ref{thmab}, we also get a solution to the
strong interpolation problem.

\bt
\label{thmc6}
Let the notation be as above.
Then there is a rational function $\phi$ of bidegree at most $(\pi^1+\nu^1 + \d^1, \pi^2+\nu^2 + \d^2)$
that solves the strong interpolation
problem  and  is unimodular on $\T^2$ except for at most a finite set.
\et

Remark: If we knew that $\phi$ could be written as the ratio of two rational inner functions,
$\phi = f/g$, then the bidegree of $f$ is at most $(\pi^1+\d^1, \pi^2 + \delta^2)$
and the bidegree of $g$ is at most $(\nu^1 + \d^1, \nu^2 + \d^2)$,
just by counting zeroes on distinguished varieties, as in \cite{agmc_dv}.

\begin{question}
When can an interpolation problem on the bidisk be solved by a ratio of rational inner functions?
\end{question}

\bibliography{../references}

\end{document}